\documentclass[12pt,a4paper,titlepage]{article} 
\usepackage{latexsym}
\usepackage{amssymb}
\usepackage{amsmath}

\def\N{{\mathbb N}}

\def\pf{\emph{Proof. }}
\def\qed{$\blacksquare$}

\def\e{\varepsilon}

\newcommand{\clspan}{\ensuremath{\overline{\operatorname{span}}}}

\def\la{\langle}
\def\ra{\rangle}
\def\be{\begin{equation}}
\def\ee{\end{equation}}

\newtheorem{thm}{Theorem}[section]
\newtheorem{lem}[thm]{Lemma} 
\newtheorem{cor}[thm]{Corollary}

\newtheorem{defn}[thm]{Definition}

\numberwithin{equation}{section}

\begin{document}

\title{Characterizing Arbitrarily Slow Convergence\\ 
in the Method of Alternating Projections}\author{Heinz H. Bauschke\footnote{Supported in part by the Natural Sciences and Engineering Research Council of Canada and by the Canada Research Chair Program.}, Frank Deutsch, and Hein Hundal }
\date{\today\; Version}
\maketitle

\begin{abstract}
Bauschke, Borwein, and Lewis\ have stated a trichotomy theorem  \cite[Theorem 5.7.16]{bbl;97} that characterizes when the convergence of the method of alternating projections can be arbitrarily slow. However,  there are  two errors in their proof of this theorem. In this note, we show that although one of the errors is critical,  the theorem itself is correct. We give a  different proof that uses the multiplicative form of the spectral theorem, and the theorem holds in any real or complex Hilbert space, not just in a real Hilbert space.

\vspace{2.5in}
\end{abstract}

\section{Introduction}\label{S: intro}

For the notation and basic Hilbert space results necessary to read this paper, the book  \cite{deu;01} is a good source, especially chapter 9.

Let $H$ be a (real or complex) Hilbert space with inner product $\la x, y \ra$ and norm $\|x\|=\sqrt{\la x, x\ra}$. If $M$ is any closed (linear) subspace of $H$, let $P_M$ denote the  orthogonal projection onto $M$. That is, $P_M: H \to M$ is defined by 
\[
\|x-P_M(x)\|=\inf_{y\in M}\|x-y\|.
\]

Let $M_1$ and $M_2$ be closed subspaces in $H$ and $M:=M_1\cap M_2$. It is well-known that $P_{M_1}P_{M_2}=P_M$ if and only if $P_{M_1}$ and $P_{M_2}$ commute: $P_{M_1}P_{M_2}=P_{M_2}P_{M_1}$. Von Neumann established the following result  which yields an interesting analogue  in the non-commuting case.

\begin{thm}\label{vNH}{\rm\textbf{(von Neumann \cite{vN;50})}} For each $x\in H$, there holds
\begin{equation}\label{vNH;eq1}
\lim_{n\to \infty} \|(P_{M_2}P_{M_{1}})^n(x)-P_M(x)\|=0.
\end{equation}
\end{thm}

The method of constructing the sequence $(P_{M_2}P_{M_{1}})^n(x)$ by alternately projecting onto one subspace and then the other is called the \emph{method of alternating projections}. While Von Neumann's theorem shows that the sequence of iterates $(P_{M_2}P_{M_{1}})^n(x)$,  \emph{always} converges to $P_M(x)$ for every $x$, it does not say anything about the speed or rate of convergence. To say something about  this, we will use the notion of angle between subpaces. Recall that the (Friedrichs) \textbf{angle} between the subspaces $M_1$ and $M_2$ is defined to be the angle  in $[0, \pi/2]$ whose cosine is given by
\[
c(M_1, M_2):=\sup\{|\la x, y\ra | \mid x\in M_1\cap M^\perp\cap B_H, \; y\in M_2\cap M^\perp\cap B_H \},
\]
where  $B_H:=\{x\in H \mid \|x\| \le 1\}$ is the unit ball in $H$.  It is  easy to see that  $0\le c(M_1, M_2) \le 1$.

\begin{thm}\label{AKW}{\rm\textbf{(Aronszajn \cite{aro;50})}} For each $x\in H$ and $n\ge 1$, we have
\begin{equation}\label{AKW;eq1}
\|(P_{M_2}P_{M_1})^n(x)-P_M(x)\| \le c(M_1, M_2)^{2n-1}\|x\|.
\end{equation}
\end{thm}

 Kayalar and Weinert \cite{kw;88} showed that the constant in Aronszajn's theorem is smallest possible independent of $x$. More precisely, they proved that 
 \begin{equation}\label{kw:eq1} 
 \|(P_{M_2}P_{M_1})^n-P_M\|=c(M_1, M_2)^{2n-1} \mbox{\quad for each $n\in \N$. }
 \end{equation}

The usefulness of the bound in (\ref{AKW;eq1}) depends on knowing when the cosine of  the angle between $M_1$ and $M_2$ is less than one, i.e., when the angle is positive.  A useful characterization of when this happens is the following.

\begin{lem}\label{DS} $c(M_1, M_2)<1$ if and only if $M_1+M_2$ is closed.
\end{lem}

This lemma is a consequence of results of Deutsch \cite{deu;84} and Simonic, whose result appeared in \cite[Lemma 4.10]{bb;93} (see also \cite[Theorem 9.35, p. 222]{deu;01}).

Recall  that a sequence $(x_n)$ is said to  converge to $x$ \textbf{linearly} provided there exists an $\alpha<1$ and a constant $c$  such that 
\[
\|x_n-x\| \le c\alpha^n \mbox{\quad for each $n\ge 1$}.
\]
In this case, we say that the rate of convergence is $\alpha$.

Using Lemma \ref{DS} and Theorem \ref{AKW}, we see that there is \emph{linear convergence} for the method of alternating projections whenever the sum of the subspaces is closed. What can be said when the sum is not closed?

Franchetti and Light \cite{fl;86} gave the first example of a Hilbert space  and two closed subspaces whose sum was not closed such that: given any sequence of reals decreasing to zero, there exists a point in the space with the property that the convergence in the von Neumann theorem was at least as slow as this sequence of reals. But this still left open the question of whether such a construction could be made in \emph{any} Hilbert space whenever $M_1$ and $M_2$ were \emph{any} closed subspaces whose sum was not closed.

In their study of the method of alternating projections, Bauschke, Borwein, and Lewis \cite{bbl;97} stated the following dichotomy. (Actually, they stated their result as a trichotomy since they were considering the more general setting of closed \emph{affine} sets, i.e., translates of subspaces,  rather than subspaces. In this situation, unlike the subspace case, one must also consider the possibility  that the intersection of the affine sets is empty. However, when the intersection is nonempty, the affine sets case easily reduces to the subspace case by a simple translation.) Roughly speaking, it states that in the method of alternating projections, either there is linear convergence for each starting point, or there exists a point which converges arbitrarily slowly.

\begin{thm}\label{BBL} {\rm (dichotomy)} Let $M_1$ and $M_2$ be closed subspaces in a  Hilbert space $H$ and $M=M_1\cap M_2$. Then exactly one of the following alternatives holds.
\begin{enumerate}
\item[{\rm(1)}] $M_1+M_2$ is closed. Then for each $x\in H$, the sequence $(P_{M_2}P_{M_1})^n(x)$ converges linearly to $P_M(x)$ with a rate  $[c(M_1, M_2)]^2$.
\item[{\rm(2)}] $M_1+M_2$ is not closed. Then for each $x\in H$, the sequence $(P_{M_2}P_{M_1})^n(x)$ converges to $P_M(x)$. But convergence is ``arbitrarily slow'' in the following sense: for each sequence $(\lambda_n)$ of positive real numbers with $1>\lambda_1 \ge \lambda_2 \ge \cdots \ge \lambda_n \to 0$, there exists a point $x_\lambda \in H$ such that 
\[
\|(P_{M_2}P_{M_1})^n(x_\lambda)-P_M(x_\lambda)\| \ge \lambda_n\mbox{\quad for all $n$.}
\]
\end{enumerate}
\end{thm}

\noindent \textbf{Remark}  Clearly, the first statement of Theorem \ref{BBL} is an immediate consequence of Theorem \ref{AKW} and Lemma \ref{DS}.  Thus we need only verify the second statement. We will do this in Section \ref{S: BBL} below.

\section{Multiplicative form of the spectral theorem}

The main fact that we will use in the proof of Theorem \ref{BBL}  is the multiplicative form of the spectral theorem (see Halmos \cite{hal;63} or Reed-Simon \cite[Corollary on p. 227]{rs;72}). Recall that a bounded linear operator $U:H_1\to H_2$ between Hilbert spaces $H_1$ and $H_2$ is called \emph{unitary} if $U$ is invertible and  $U^*=U^{-1}$. It follows that a unitary operator  is isometric: $\|Ux\|=\|x\|$ for each $x\in H$. Since the inverse of a unitary operator is unitary, it too is isometric. (We will use these facts in a few places below without explicit mention.)

\begin{thm}\label{RS}{\rm\textbf {(Spectral Theorem; multiplicative form)}} Let $H$ be a  (real or complex) Hilbert space, and let $T$ be a self-adjoint bounded  linear operator on $H$. Then there exists a finite measure space $(\Omega, \mu )$, a bounded real-valued function $F$ on $\Omega$, and a unitary map $U:  H\to L_2(\Omega, \mu)$ such that 
\begin{equation}\label{RS;eq1}
UTU^{-1}f=F\cdot f \mbox{\quad for all $f\in L_2(\Omega, \mu)$.}
\end{equation} 
Defining  $D: L_2(\Omega, \mu) \to L_2(\Omega, \mu)$ to be the operator ``multiplication by $F$'', $(Df)(t):=F(t)f(t)$,  this can be expressed in operator notation as 
\begin{equation}\label{RS;eq2}
UTU^{-1}=D.
\end{equation}
\end{thm}

 Actually, in both \cite{hal;63} and \cite{rs;72}, the theorem is stated for a \emph{complex} Hilbert space only, and \cite{rs;72} even assumes separability. However,  it is easy to check that each of the tools used in the proof in \cite{hal;63}, for example, has a corresponding real space analogue.

\medskip\noindent
\textbf{Acknowledgements} We  are greatly indebted to Joel Anderson, Nigel Higson, and Barry Simon for personally transmitting some very useful comments to us related to the multiplicative form of the spectral theorem.

A self-adjoint operator $T$ on $H$ is called \emph{positive} if $\la Tx, x\ra \ge 0$ for each $x\in H$. A simple, but important, example of a positive operator is the orthogonal projection $P_S$ onto any closed subspace $S \subset H$ (see, e.g., \cite[p. 79]{deu;01}).

\begin{cor}\label{cor to RS} Assume the hypothesis of Theorem \ref{RS}. If $T$ is also positive, then the bounded real-valued function $F$ of Theorem \ref{RS} is also nonnegative a.e.$(\mu)$.
\end{cor}

\pf  Let $f\in L_2(\Omega, \mu)$ be arbitrary and $y=U^{-1}f$. Since $T$ is positive, we have that
\begin{eqnarray*}
\int _{\Omega}F|f|^2d\mu&=& \la Ff, f\ra =\la Df, f\ra= \la UTU^{-1}f, f\ra\\
&=&\la TU^{-1}f, U^*f\ra =\la Ty, y\ra \ge 0.
\end{eqnarray*} 
Briefly, $\int _{\Omega}F|f|^2d\mu \ge 0$ for each $f\in  L_2(\Omega, \mu)$. We readily deduce that $F\ge 0$ a.e.($\mu$).  \qed

\section{Proof of Theorem \ref{BBL} }\label{S: BBL}

In this section we will prove the second statement of Theorem \ref{BBL}.    Our proof is along the same general  lines as in  \cite{bbl;97} in that we proceed by a series of small  steps that are each easily digested.  However, there are subtle errors in steps 2 and 3 of \cite{bbl;97} (see Section  \ref{S:Errors} for the details). We will avoid these errors by using Theorem \ref{RS} and following a somewhat  different path.

\textbf{Proof of the second statement in Theorem \ref{BBL}.}  Suppose  $M_1+M_2$ is not closed, and let $(\lambda_n)$ be a sequence with $1>\lambda_1\ge \lambda_2 \ge \cdots \ge \lambda_n>0$, and $\lambda_n \to 0$.  By Lemma \ref{DS}, $c(M_1, M_2)=1$. Let
\begin{equation}\label{eq1}
A=M_1\cap M^\perp \mbox{\quad and \quad} B=M_2\cap M^\perp.
\end{equation}
Note that $A$ and $B$ are closed subspaces with $A\cap B=\{0\}$. Clearly,
\be\label{eq1}
c(A, B)=c(M_1, M_2)=1
\ee
and hence, by Lemma \ref{DS} again, $A+B$ is not closed. Since $c(A,  B)=\|P_BP_A\|$ by \cite{deu;84} (see also \cite[Lemma 9.5(7), p. 197]{deu;01}), it follows that $\|P_BP_A\|=1$.

\begin{lem} \label{T facts} The operator $T:=P_AP_BP_A$ is a bounded self-adjoint linear operator on $H$ which is positive and $\|T\|=1$. Hence there exists a finite measure space $(\Omega, \mu)$, a nonnegative bounded function $F$ on $\Omega$, and a unitary operator $U: H\to L_2:=L_2(\Omega, \mu)$ such that
\begin{equation}\label{eq1}
UTU^{-1}=D,
\end{equation}
where $D: L_2 \to L_2$ is defined by $Df:=Ff$ for each $f\in L_2$.
\end{lem}

\emph{Proof of Lemma \ref{T facts}}. By Corollary \ref{cor to RS}, it suffices to verify the first statement of the lemma. Clearly, $T$ is self-adjoint and bounded.  Moreover, using \cite[Corollary 5.17]{dhII;06}, $\|T\|=\|P_AP_BP_A\|=\|P_BP_A\|^2=1$.  Fix any $x\in H$ and set  $y=P_Ax$. Since $P_B$ is positive, we have that
\[
\la Tx, x \ra=\la P_AP_BP_Ax, x \ra= \la P_BP_Ax, P_Ax \ra= \la P_By, y\ra \ge 0.
\]
 This shows that $T$ is positive on $H$ and completes the proof of Lemma \ref{T facts}.
\medskip

For each $k\in \N\!:=\{1, 2, \dots\}$, let $s_k$ be the largest integer such that $s_k\lambda_k<1$. Then the following claim is clear.

\medskip
\textbf{Claim 1.}  \textsl{$s_k\lambda_k<1 \le (s_k+1)\lambda_k$ for all $k\in \N$, $s_1\le s_2\le s_3\le \cdots$, and each $s_k$ occurs only finitely often.}
\medskip

Next let $(t_n)$ be the strictly increasing sequence of integers with 
\begin{equation}
\{t_1, t_2, \dots \}=\{s_1, s_2, \dots \}.
\end{equation}

Note that since $(t_n)$ is a subsequence of $(n)$, it follows that 
\be
\sum_1^\infty \frac1{t_n^2} < \infty.
\ee

For each $n\in\N$, we define
\be
k_0(n): =\min \{ k\mid s_k=t_n \} \mbox{\quad  and\quad } k_1(n):=\max \{ k \mid s_k=t_n\}.
\ee

It is clear that $k_0(n) \to \infty$, $k_1(n) \to \infty$, and
\be\label{eq1;claim 1}
s_{k_0(n)-1}=t_{n-1}<t_n=s_{k_0(n)}=s_{k_0(n)+1}=\cdots=s_{k_1(n)} < t_{n+1}=s_{k_1(n)+1}.
\ee

Set 
\be\label{eq2;claim 1}
\alpha_n:=(\lambda_{k_0(n)}t_n)^{\frac1{2k_1(n)}} \mbox{\quad for each $n\in \N$}.
\ee

\medskip
\textbf{Claim 2.} \textsl{ For each $n\in \N$,
\be\label{eq1;claim 2}
1> \lambda_{k_0(n)}s_{k_0(n)}=\lambda_{k_0(n)}t_n \ge 1-\lambda_{k_0(n)},
\ee
\be\label{eq2;claim 2}
0< \alpha_n < 1, \mbox{\quad and \quad} \alpha_n \to 1.
\ee}

To see this, note that by definition, $\lambda_{k_0(n)}t_n=\lambda_{k_0(n)}s_{k_0(n)} < 1$, and $1\le \lambda_{k_0(n)}(s_{k_0(n)}+1)$. But the latter inequality implies that $1-\lambda_{k_0(n)} \le \lambda_{k_0(n)}s_{k_0(n)}=\lambda_{k_0(n)}t_n$. Also, $\lambda_{k_0(n)}t_n < 1$ implies that $\alpha_n < 1$. Since $\lambda_{k_0(n)} \to 0$, relation (\ref{eq1;claim 2}) implies that $\lambda_{k_0(n)}t_n \to 1$. This, along with $k_1(n) \to \infty$, shows that $\alpha_n \to 1$, which completes the proof of Claim 2.

We note that the first two claims follow exactly as in the proof given in \cite{bbl;97}. However, at this point our approach will deviate significantly from that of \cite{bbl;97}.

\medskip
\textbf{Claim 3.}  $\mu\{F^{-1}([1, \infty))\} =0$.

\medskip
To see this, let $S:=F^{-1}[1, \infty)$ and $y=U^{-1}(\chi_S)$, where $\chi_S$ denotes the characteristic function of $S$: $\chi_S(t)=1$ if $t\in S$ and $0$ otherwise. We must show that $\mu(S)=0$. Since 
\be\label{first eq; claim 3}
\|y\|= \|U^{-1}(\chi_S)\|=\|\chi_S\|=\left(\int_S1d\mu\right)^{1/2}=[\mu(S)]^{1/2},
\ee
it suffices to show that $y=0$. Using (\ref{first eq; claim 3}), we have
\begin{eqnarray}
\|Ty\|&=& \|U^{-1}DUy\|=\|U^{-1}D( \chi_{S})\| = \|U^{-1}(F\chi_S)\|=\|F\chi_S\| \nonumber \\
&=&\left[\int_SF^2d\mu\right]^\frac12 \ge \left[\int_S1d\mu\right]^\frac12=\|y\|. \label{eq1;claim 3}
\end{eqnarray}

This shows  that $\|Ty\|\ge \|y\|$. But since $T=P_AP_BP_A$ is the product of norm one operators, $\|Ty\|\le \|y\|$. Thus $\|Ty\|=\|y\|$. We deduce that
\be\label{eq2;claim 3}
\|y\|=\|P_AP_BP_Ay\| \le \|P_BP_Ay\| \le \|P_Ay\| \le \|y\|.
\ee
Thus we must have equality holding throughout the string of  inequalities (\ref{eq2;claim 3}). It follows (see, e.g., \cite[Theorem 5.8(2), p. 76]{deu;01}) that $y\in A\cap B=\{0\}$ and hence $y=0$. This proves Claim 3.

\medskip
\textbf{Claim 4.} \textsl{For each} $\e>0$, $\mu\{F^{-1}((1-\e, 1))\}>0$.

\medskip
If not, there exists $\e>0$ such that  $\mu\{F^{-1}((1-\e, 1))\}=0$. Choose any $ y \in H$ and set $g=Uy$. Then, using Claim 3, we have that
\begin{eqnarray*}
\|Ty\|^2&=&\|U^{-1}DUy\|^2=\|DUy\|^2=\|Dg\|^2=\int|Fg|^2d\mu =\int F^2|g|^2 d\mu\\
&=&\int_{F^{-1}([0, 1-\e])}F^2|g|^2d\mu +\int_{F^{-1}((1-\e, 1))}F^2|g|^2d\mu +\int_{F^{-1}([1, \infty))}F^2|g|^2d\mu \\
&\le&(1-\e)^2\int_{F^{-1}([0, 1-\e])}|g|^2d\mu +0 +0 \le (1-\e)^2\int |g|^2d\mu\\
&=& (1-\e)^2\|g\|^2=(1-\e)^2\|Uy\|^2=(1-\e)^2\|y\|^2.
\end{eqnarray*}
Briefly, $\|Ty\| \le (1-\e)\|y\|$ for each $y\in H$. It follows that $\|T\| \le 1-\e$,  which (by Lemma \ref{T facts}) contradicts $\|T\|=1$. This proves Claim 4.

\medskip
\textbf{Claim 5.} \textsl{For each $\e>0$, there exists $\e_1 \in (0, \e)$ such that} 
\[
\mu\{ F^{-1}((1-\e, 1-\e_1))\} >0.
\]

To verify this, we use Claim 4 and the countable additivity of $\mu$ to obtain
\begin{eqnarray*}
0&<& \mu\{F^{-1}((1-\e, 1))\}=\mu\left\{\bigcup_{i=1}^\infty F^{-1}\left(\left(1-\frac{\e}i, 1-\frac{\e}{i+1}\right]\right)\right\}\\
&=&\sum_{i=1}^\infty \mu\left\{F^{-1}\left(\left(1-\frac\e{i}, 1-\frac\e{i+1}\right]\right)\right\}.
\end{eqnarray*}
Thus there exists an integer $i$ such that $\mu\left\{F^{-1}\left((1-\frac{\e}{i}, 1-\frac{\e}{i+1}]\right)\right\} >0$. Let $\e_1=\frac{\e}{i+2}$. Then $\e_1 \in (0, \e)$ and
\[
\mu\left\{ F^{-1}\left((1-\e, 1-\e_1)\right)\right\} \ge \mu\left\{F^{-1}\left((1-\frac\e{i}, 1-\e_1)\right)\right\}  > 0.
\]
This proves Claim 5.

\medskip
\textbf{Claim 6.} \textsl{There exists a sequence  of reals $(\beta_n)\subset (0, 1)$ such that} $\alpha_n^2\le \beta_n< \beta_{n+1}<1$ and $\mu\{F^{-1}\left([\beta_n, \beta_{n+1})\right)\}>0$ \textsl{for each $n\in \N$}.

\medskip
We prove Claim 6 by induction. For $n=1$, take $\beta_1=\alpha_1^2$. Then $\beta_1 < 1$. Assume next that $\beta_1, \dots, \beta_m$ have been chosen so that $\beta_1<\beta_2< \cdots <\beta_m < 1$, $\beta_k \ge \alpha_k^2$ for $k=1, 2, \dots, m$, and $\mu\{F^{-1}\left( [\beta_k, \beta_{k+1})\right)> 0$ for $k=1, 2, \dots, m-1$. Let $\e:=\min\{ 1-\alpha^2_{m+1}, 1-\beta_m\}$. Then $\e>0$ and Claim 5 implies the existence of $\e_1\in (0, \e)$ such that $\mu\{F^{-1}\left((1-\e, 1-\e_1)\right)\}>0$.  Let $\beta_{m+1}:=1-\e_1$. Then $\beta_{m+1}> 1-\e \ge \beta_m$. Also, $\beta_{m+1} > 1-\e \ge \alpha^2_{m+1}$. Finally, $\mu\{F^{-1}\left( [\beta_m, \beta_{m+1})\right) \} \ge \mu\{F^{-1}\left( [1-\e, 1-\e_1)\right)\} >0$. This completes the induction step and hence the proof.

\begin{defn} \label{orthog vectors} With  $\beta_n$ given as in Claim 6, for each $n\in \N$,  let $S_n:=F^{-1}( [\beta_n, \beta_{n+1}))$ and define the vector $e_n\in H$ by
\[
e_n:=\frac1{\sqrt{\mu(S_n)}}U^{-1}(\chi_{S_n}).
\]
\end{defn}

Note that $$Ue_n=\frac1{\sqrt{\mu(S_n)}}\chi _{S_n}.$$

\textbf{Claim 7.} $\|e_n\|=1$  \textsl{for each} $n\in \N$.

\medskip This follows from
\[ 
\|e_n\|=\|Ue_n\|=\frac1{\sqrt{\mu(S_n)}}\|\chi_{S_n}\|=1.
\]
\medskip

It is convenient to list next a few basic and easily verified facts concerning powers of $T$ and $D$.

\medskip \textbf{Claim 8.}\textsl{\begin{enumerate}
\item[{\rm (1)}] $T^k=(U^{-1}DU)^k=U^{-1}D^kU$.
\item[{\rm (2)}]  $D^kf=F^kf$ for all $f\in L_2(\Omega, \mu)$.
\item[{\rm (3)}]  If $f, g \in L_2(\Omega, \mu)$ and $f(t)g(t)=0$ for $\mu$ almost all $t$, then $\la D^jf, D^kg\ra =0$ for every $ j, k \in \N\cup \{0\}$.
\end{enumerate}}

\medskip
\textbf{Claim 9.} \textsl{For all integers $j, k \in \N\cup\{0\}$ and $m, n \in \N$ with $m\ne n$, we have
\[
\la T^je_m, T^ke_n\ra =0.
\] }

\medskip
To verify this, let $f_r:=Ue_r=\frac1{\sqrt{\mu (S_r)}}\chi_{S_r}$ for each $r\in \N$. Then
$
\chi_{S_n}\chi_{S_m}=\chi_{S_n\cap S_m}=0
$
since $S_n\cap S_m=\emptyset$. Thus $f_nf_m=0$.  Using statements (1) and (3) of Claim 8, we get that
\[
\la T^je_m,  T^ke_n\ra=\la U^{-1}D^jUe_m, U^{-1}D^kUe_n\ra=\la D^jf_m, D^kf_n\ra=0.
\]

\medskip
\textbf{Claim 10.} $\beta_{n+1}^k \ge \|T^ke_n\| \ge \beta_n^k\ge \alpha_n^{2k}$ \textsl{for all} $k, n \in \N$.

\medskip

The last inequality follows from Claim 6. Next observe that 
\begin{eqnarray*}
\|T^ke_n\|^2&=&\|U^{-1}D^kUe_n\|^2=\left\|D^k(\frac{\chi_{S_n}}{\sqrt{\mu (S_n)}})\right\|^2=\int \left[D^k(\frac{\chi_{S_n}}{\sqrt{\mu (S_n)}})\right]^2d\mu \\
&=&\frac1{\mu(S_n)}\int F^{2k}\chi_{S_n}^2d\mu=\frac1{\mu(S_n)}\int_{S_n}F^{2k}d\mu.
\end{eqnarray*}
Also, by the definition of $S_n$ (in Definition \ref{orthog vectors}), it is clear that
\[
\beta_n^{2k} \le \frac1{\mu(S_n)}\int_{S_n}F^{2k}d\mu \le \beta_{n+1}^{2k}.
\]
Taking square roots completes the proof of Claim 10.

\medskip
Now we can define the element which will converge slower than the sequence $(\lambda_n)$.
\begin{defn}\label{the slow guy} Set 
\[
x_\lambda:=\sum_1^\infty\frac1{t_n} e_n.
\]
\end{defn}

 Since $\sum_1^\infty 1/t^2_n \le \sum_1^\infty 1/n^2 < \infty$ and $\|e_n\|=1$, it follows that $x_\lambda$ is a well-defined element of $H$.

\medskip
\textbf{Claim 11.} $\|T^kx_\lambda\| \ge \alpha_n^{2k}/t_n$ \textsl{for all  $n, k\in \N$.}

\medskip
We deduce
\begin{eqnarray*}
\|T^kx_\lambda\|^2&=&\la T^kx_\lambda, T^kx_\lambda\ra=\left\la T^k\left(\sum_ne_n/t_n\right), T^k\left(\sum_me_m/t_m\right) \right\ra \\
&=&\sum_n\frac 1{t_n} \sum_m\frac1{t_m} \la T^ke_n, T^ke_m\ra\\
&=&\sum_n\frac1{t_n^2}\|T^ke_n\|^2 \mbox{\quad (by Claim 9)}\\
&\ge& \frac1{t_n^2}\|T^ke_n\|^2 \mbox{\quad for each $n$}\\
&\ge & \frac{\alpha_n^{4k}}{t^2_n} \mbox{\quad (by Claim 10)}.
\end{eqnarray*}
Thus $\|T^kx_\lambda\| \ge {\alpha_n^{2k}}/{t_n}$ as claimed.

\medskip
\textbf{Claim 12.} $\|(P_BP_A)^kx_\lambda\| \ge \lambda_k \mbox{\quad \textsl{for each} $k\in \N$.}$

\medskip
Fix any $k\in \N$ and choose $n\in \N$ such that  $k_0(n) \le k \le k_1(n)$. Using Claim 11, we get that
\begin{equation}
\|(P_BP_A)^kx_\lambda\| \ge \|P_A(P_BP_A)^kx_\lambda\|=\|T^kx_\lambda\|\ge \frac{\alpha_n^{2k}}{t_n} \ge \frac{\alpha_n^{2k_1(n)}}{t_n}=\lambda_{k_0(n)}\ge \lambda_k, \nonumber
\end{equation}
which proves Claim 12.

\medskip
\textbf{Claim 13.} \textsl{For each $k\in \N$}, $(P_{M_2}P_{M_1})^k-P_M=(P_BP_A)^k$.

\medskip
Using the facts that $M=M_1\cap M_2$, $P_{M^\perp}=I-P_M$, and $P_{M^\perp}$ is idempotent and commutes with both $P_{M_1}$ and $P_{M_2}$ (see, e.g., \cite[p. 194]{deu;01}), we get that $P_{M_i}P_{M^\perp}=P_{M_i\cap M^\perp}$ for $i=1, 2$ and
\begin{eqnarray*}
(P_{M_2}P_{M_1})^k-P_M&=&(P_{M_2}P_{M_1})^k(I-P_M)=(P_{M_2}P_{M_1})^kP_{M^\perp}\\
&=&(P_{M_2}P_{M^\perp}P_{M_1}P_{M^\perp})^k=(P_{M_2\cap M^\perp}P_{M_1\cap M^\perp})^k\\
&=&(P_BP_A)^k,
\end{eqnarray*}
which proves Claim 13. 

Combining Claims 12 and 13, we immediately obtain

\medskip

\textbf{Claim 14.} $\|(P_{M_2}P_{M_1})^k(x_\lambda)-P_M(x_\lambda)\| \ge \lambda_k $\ \textsl{ for each} $k\in \N$. 

\medskip

This completes the proof of the second statement of Theorem \ref{BBL}.

\section{Two errors in \cite{bbl;97}}\label{S:Errors}

In this section, we point out two errors in \cite{bbl;97}. We shall use the notation of \cite{bbl;97}. (Note that this is the same as the notation of the present paper except that here we have used $M_1, M_2$ instead of $C_1, C_2$.)

\medskip

\textbf{First error.}
The proof of the Claim in Step~2 of the proof of Theorem~5.7.16 in \cite{bbl;97}
has a mistake. The Claim itself is correct, only the proof of this claim
is incorrect. 

Specifically, we inductively construct  $(e_n')$ and $(f_n')$ in $A$ and $B$,
respectively. 
Let $E$ and $F$ be the finite-dimensional spaces as in the proof.
Let $(a_n)$ in $A$ and $(b_n)$ in $B$ as in the proof: 
\begin{equation} \label{eq:eins}
\|a_n\|=1=\|b_n\|\quad\text{and}\quad  \la a_n, b_n\ra \to 1,
\end{equation}
and $a_n\to 0$ weakly  and $b_n\to 0$ weakly. 
Because $E+F$ is \emph{finite-dimensional}, the sum
$A^\bot + (E+F)$ is closed. Hence $\{A^\bot,E+F\}$ is regular
(by \cite[Proposition~5.16]{bb;96}) and so is $\{A^{\bot\bot},(E+F)^\bot\} =
\{A,E^\bot \cap F^\bot\}$ (again by \cite[Proposition~5.16]{bb;96}). 
This means the following by definition of regularity.

\textbf{Observation.}\textsl{
If $(z_n)$ is a bounded sequence with \\
$\max\big\{d(z_n,A), d(z_n,E^\bot\cap F^\bot)\big\}\to 0$,
then $d(z_n,A \cap E^\bot\cap F^\bot)\to 0$. 
(And analogously when $A$ is replaced by $B$.)}

Now back to the proof of the Claim. This time, $P_{E+F}$ is a compact
operator. (In \cite{bbl;97}, $P_E$ and $P_F$ were considered, which is not
sufficient.)  Since $a_n \to 0$ weakly and $b_n\to 0$ weakly, we deduce that
\begin{equation}
P_{E+F}a_n \to 0\quad\text{and}\quad P_{E+F}b_n \to 0.
\end{equation}
Since $(E+F)^\bot = E^\bot \cap F^\bot$, this implies 
\begin{equation}
a_n-P_{E^\bot \cap F^\bot}a_n \to 0 \quad\text{and}\quad
b_n-P_{E^\bot \cap F^\bot}b_n \to 0;
\end{equation}
equivalently,
\begin{equation}
d(a_n,E^\bot \cap F^\bot)\to 0 \quad\text{and}\quad
d(b_n,E^\bot \cap F^\bot)\to 0.
\end{equation}
The above Observation now implies 
$d(a_n,A\cap E^\bot \cap F^\bot)\to 0$ and 
$d(b_n,B\cap E^\bot \cap F^\bot)\to 0$; equivalently,
\begin{equation}
a_n-P_{A \cap E^\bot \cap F^\bot}a_n \to 0 \quad\text{and}\quad
b_n-P_{B \cap E^\bot \cap F^\bot}b_n \to 0.
\end{equation}
In view of \eqref{eq:eins}, we deduce that
\begin{equation}
\la P_{A \cap E^\bot \cap F^\bot}a_n, P_{B \cap E^\bot \cap F^\bot}b_n\ra
\to 1.
\end{equation}
Thus, for all $n$ sufficiently large,
we have $\|P_{A \cap E^\bot \cap F^\bot}a_n\|\leq 1$,
$\|P_{B \cap E^\bot \cap F^\bot}b_n\|\leq 1$,
$P_{A \cap E^\bot \cap F^\bot}a_n\in A \cap E^\bot \cap F^\bot$,  
$P_{B \cap E^\bot \cap F^\bot}b_n\in B \cap E^\bot \cap F^\bot$,
and \\
$\la P_{A \cap E^\bot \cap F^\bot}a_n, P_{B \cap E^\bot \cap
F^\bot}b_n \ra$ is as close to $1$ (from below) as we like. Then for $n$ sufficiently large, we can take $e_{m+1}'=P_{A\cap E^\perp\cap F^\perp}a_n$ and $f_{m+1}'=P_{B\cap E^\perp\cap F^\perp}b_n$.

\medskip
\textbf{Second error.}
The second error is on the third line on page 32 of \cite{bbl;97}, where it is claimed that
\begin{equation} \label{e:falsch}
C_1 = (C_1\cap C_2)\oplus E \oplus (A \cap E^\bot \cap F^\bot),\;
C_2 = (C_1\cap C_2)\oplus F \oplus (B \cap E^\bot \cap F^\bot).
\end{equation}
 Unfortunately, only 
\begin{equation*}
C_1 = (C_1\cap C_2)\oplus E \oplus (A\cap E^\bot),\quad
C_2 = (C_1\cap C_2)\oplus F \oplus (B\cap F^\bot)
\end{equation*}
is true. This invalidates the rest of the proof in \cite{bbl;97}.

Here is a counterexample to \eqref{e:falsch}.
Let $\{u_n \mid n\in \N\}$ be an orthonormal basis of a separable Hilbert space. 
Set 
\begin{equation}\nonumber
C_1 := \clspan \{ u_{2n}+\tfrac{1}{n}u_{2n-1}\mid n\in \N\}
\quad\text{and}\quad
C_2 := \clspan \{ u_{2n}+\tfrac{1}{n}u_{2n+1}\mid n\in \N \}.
\end{equation}
Then 
\begin{equation}
C_1 \cap C_2 =\{0\}. 
\end{equation}
(Sketch: the spanning vectors are orthogonal. Normalize and use Fourier
expansions. Equate coefficients, compare odd and even ones. Deduce that
they are all equal; thus they must be equal to $0$.)
Hence $A=C_1$ and $B=C_2$. 
Set
\begin{equation}
e_n = e_n' = \rho_n\big(u_{4n} + \tfrac{1}{2n}u_{4n-1}\big)\quad\text{and}\quad
f_n = f_n' =\rho_n \big(u_{4n} + \tfrac{1}{2n}u_{4n+1}\big),
\end{equation}
where $\rho_n:=(1+\tfrac{1}{4n^2})^{-1/2}$.
Since $\la {e_n'}, {f_n'} \ra= (1+\tfrac{1}{4n^2})^{-1}$,
the sequences $(e_n')$ and $(f_n')$ are as in the Claim of Step~2,
and the sequences $(e_n)$ and $(f_n)$ are as in Step~3. 
Set 
\begin{equation}
E = \clspan\{e_n \mid n\in \N\}
\quad\text{and}\quad
F = \clspan\{f_n \mid n\in\N\}.
\end{equation}
Then
\begin{equation}
\overline{E+F} = \clspan\{2ne_{4n}+u_{4n-1}, \; 2ne_{4n}+u_{4n+1} \mid n\in \N \}
\end{equation}
is a subspace of $\clspan\{u_{4n-1},u_{4n},u_{4n+1}\mid n\in \N \}$.
Thus $\{u_1,u_2,u_6,u_{10},\ldots\}\subset (E+F)^\bot$. 
Since the orthogonal complement of 
$\clspan\{2ne_{4n}+u_{4n-1},2ne_{4n}+u_{4n+1} \mid n\in \N\}$ in
$\clspan\{u_{4n-1},u_{4n},u_{4n+1} \mid n\in \N\}$ is
$\clspan\{-2nu_{4n-1}+u_{4n} - 2nu_{4n+1}\mid n\in \N \}$, we obtain 
\begin{eqnarray} \label{e:EbotFbot}
E^\bot \cap F^\bot &=&(E+F)^\bot  \\
&=& \clspan\{u_1,\; u_{4n-2},\; -2nu_{4n-1}+u_{4n} -
2nu_{4n+1}\mid n\in \N \}. \nonumber
\end{eqnarray}
Consider the vector $x := u_6 + \tfrac{1}{3}u_5$.
Then $x$ belongs to $C_1 = A$. Since $E \subset
\clspan\{u_{4n-1},u_{4n} \mid n\in \N \}$, it follows that $x\in E^\bot$
and hence $P_Ex = 0$. 
Now consider the first term in the false statement \eqref{e:falsch}, which
in our present situation becomes
\begin{equation} \label{e:falsch2}
A = E \oplus (A \cap E^\bot \cap F^\bot).
\end{equation}
This would imply that $x$ belongs entirely to $A \cap E^\bot \cap F^\bot$. 
While it is true that $x\in A \cap E^\bot$, it is \emph{not} true that
$x$ belongs to $E^\bot \cap F^\bot$. This can be verified using relation  \eqref{e:EbotFbot}.

\bibliographystyle{plain}

\begin{tabular}{lll}
Heinz H. Bauschke &Frank Deutsch    &  Hein Hundal\\
Mathematics &Department of Mathematics  &   146 Cedar Ridge Drive\\
UBC Okanagan &Penn State University  & Port Matilda, PA 16870\\
Kelowna, British Columbia &University Park, PA 16802  & {} \\
V1V 1V7, Canada & USA &USA\\
 heinz.bauschke@ubc.ca &  deutsch@math.psu.edu &  hundalhh@yahoo.com
\end{tabular}

\end{document}